\documentclass[lineno]{biometrika}

\usepackage{amsmath}

\usepackage{times}
\usepackage{bm}
\usepackage{natbib}

\usepackage[plain,noend]{algorithm2e}

\makeatletter
\renewcommand{\algocf@captiontext}[2]{#1\algocf@typo. \AlCapFnt{}#2} % text of caption
\def\@algocf@capt@plain{top}
\renewcommand{\algocf@makecaption}[2]{%
  \addtolength{\hsize}{\algomargin}%
  \sbox\@tempboxa{\algocf@captiontext{#1}{#2}}%
  \ifdim\wd\@tempboxa >\hsize%     % if caption is longer than a line
    \hskip .5\algomargin%
    \parbox[t]{\hsize}{\algocf@captiontext{#1}{#2}}% then caption is not centered
  \else%
    \global\@minipagefalse%
    \hbox to\hsize{\box\@tempboxa}% else caption is centered
  \fi%
  \addtolength{\hsize}{-\algomargin}%
}
\makeatother

\def\Bka{{\it Biometrika}}
\def\AIC{\textsc{aic}}
\def\T{{ \mathrm{\scriptscriptstyle T} }}
\def\v{{\varepsilon}}

\usepackage[T1]{fontenc}
\usepackage{amsmath, amssymb, amsfonts}
\usepackage{dsfont}
\usepackage{natbib}
\usepackage{graphicx}
\usepackage{slashbox}
\usepackage{color}
\usepackage{aeguill}
\usepackage{multicol}
\usepackage{multirow}
\usepackage{fancyhdr}
\usepackage{pstricks}
\usepackage{enumerate}

\usepackage{bm}

\newcommand{\norm}[1]{\|#1\|}

\newcommand{\bs}{\boldsymbol}
\newcommand{\TV}{\textsc{tv}}

\newcommand{\E}{\mathbb E}

\newcommand{\argmin}{\text{argmin}}

\renewcommand{\H}{{\mathbf H}}
\newcommand{\card}{\mathrm{Card}}

\newcommand{\ind}[1]{\mathbf 1{(#1)}}%

\begin{document}

\jname{}
%% The year, volume, and number are determined on publication
%\jyear{}
%\jvol{}
%\jnum{}
%% The \doi{...} and \accessdate commands are used by the production team
%\doi{10.1093/biomet/asm023}
%\accessdate{}
%\copyrightinfo{}

%% These dates are usually set by the production team
%\received{}
%\revised{}

%% The left and right page headers are defined here:
%\markboth{}{Biometrika style}

%% Here are the title, author names and addresses
\title{A penalized algorithm for event-specific rate models for recurrent events}

\author{O. Bouaziz}
\affil{MAP5, UMR CNRS 8145 and University Paris Descartes, Paris, France, \email{olivier.bouaziz@parisdescartes.fr}}

\author{\and A. Guilloux}
\affil{LSTA, University Pierre et Marie Curie, 4 place Jussieu, Paris, France,\email{agathe.guilloux@upmc.fr}}

\maketitle

\begin{abstract}
%\texttt{ A faire}
%\begin{enumerate}
%\item OK changer l'abstract
%\item discussion 1. AFT, 2. LASSO, 3. distance finie
%\item Supp Mat.
%\item raccourcir
%\item ajouter "Further details are provided in Supplementary Material." quad il faut
%\item Agathe : Regarder les paragraphes où on dit la même chose (p2,3,4,6)
%\item Olivier : regarder les 2 lambdas. \texttt{Pour le modèle multiplicatif on ne peut pas faire mieux que le modèle constant en terme d'AIC ou BIC si on ne prend pas l'estimateur adaptatif. J'ai mis en supplementary material l'estimateur TV non adaptatif (j'ai choisi} $\lambda$ \texttt{pour avoir le meilleur AIC parmi tous les estimateurs possibles, sauf le constant). Il faut que je calcule la vraisemblance dans les modèles additifs, on ne l'a pas encore fait. Je vais voir si je m'en sors tout seul.}
%\item OK Agathe : faire les autres simulations
%\item OK Agathe : changer valeur de eps=0.001 dans calcul.adaptive.TV
%\end{enumerate}
We introduce a covariate-specific total variation penalty in two semiparametric models for the rate function of recurrent event process. The two models are a stratified Cox model, introduced in~\cite{prentice}, and a stratified Aalen's additive model. We show the consistency and asymptotic normality of our penalized estimators. We demonstrate, through a simulation study, that our estimators outperform classical estimators for small to moderate sample sizes. Finally an application to the bladder tumour data of~\cite{byar} is presented.

\end{abstract}


\section{Introduction}
%\texttt{ATTENTION TOUT AU PRESENT, VERIFIER POUR L4INTRO}

Recurrent events are frequent in clinical or epidemiological studies when each subject experiences repeated events over the time. Standard medical examples include the repetition of asthma attacks, epileptic seizures or tumour recurrences for individual patients. 
In this context, proportional hazards models have been largely studied in the literature to model the rate or mean functions of recurrent event data. For instance, \cite{Andersenetal82} introduce a conditional Cox model where the recurrent events process is assumed to be a Poisson process. Without this assumption, similar proportional hazards models and extensions are considered in~\cite{lawlessNadeau},~\cite{linetal},~\cite{linwei00} and~\cite{caischaubel}. %A stratified Event specifics models can allow more flexibility to the data through strafication

To model rate functions in a recurrent events context, a different approach consists in fitting a Cox model for any different recurrence. Along these lines,~\cite{prentice} introduce two stratified proportional hazards models with event-specifics baseline hazards and regression coefficients. 
%To this end, two stratified proportional hazards models were introduced by~\cite{prentice} with event specifics baseline hazards and regression coefficients. 
Gap times and conditional models are presented in their paper and a marginal event-specific model is studied in~\cite{weietal}. 
We refer to~\cite{kelly} for a complete review of existing Cox-based recurrent event models.% In the latter, two new variant models were also introduced: an event-specific total time model and an unrestricted gap time model. 

Additive models provide an useful alternative to proportional hazards models. For classical counting processes, the Aalen model was first introduced in~\cite{aalen} and is extensively studied in~\cite{mckeague88},~\cite{huffer91},~\cite{linying94}. It is considered in the context of recurrent events in~\cite{scheike}. We propose in this paper to consider an event-stratified version of the Aalen model, in the manner of~\cite{prentice}.
As demonstrated in the following, event-stratified models allow more flexibility but suffer from over-parametrization as soon as the sample size is not large enough with respect to the number of covariates and the number of recurrent events. We address this drawback by introducing new estimators defined as minimizers of penalized empirical risks. More specifically, we consider a covariate-specific total variation penalty. %Our choice of penalization takes advantage of the chronological nature of recurrent events and is the expression of a preference for simple models.

  The remainder of this article is organized as follows. The multiplicative and additive models studied in this paper are presented in Section~\ref{eq:modelmult}. In Paragraph~\ref{par:TVpen}, we describe our novel algorithms. It requires preliminary details on inference in these two models, which are given in Paragraphs~\ref{par:detailmult} and~\ref{par:detailadd}. Consistency and asymptotics normality of the estimators are derived in Section~\ref{sec:asymp}. Simulation studies and a real data analysis are provided in Sections~\ref{sec:simu} and~\ref{sec:real}.  A discussion and some concluding remarks are contained in Section~\ref{sec:discussion}. 
  
\section{Models and algorithm}\label{sec:model}
\subsection{Models}
Let $D$ denote the time of the terminal event and $N^*(t)$ the number of recurrent events before time $t$. The end-point of the observation is $\tau >0$.  The $p$-dimensional process of covariates is denoted by $X$ and $\rho_0$ represent the rate function. The event-specific rate function of the process $N^*$ is then defined as
\begin{equation*}
\E \big( dN^*(t)\mid X(t), D\geq  t,N^*(t)=s-1\big)=\ind{D \geq t} \rho_0(t,s,X(t)) dt,
\end{equation*}
%\texttt{Je pense qu'il faut plutot écrire
%\begin{equation*}
%\E \big( dN^*(t)\mid X(t), D\geq  t,N^*(t)=s-1\big)=\ind{D \geq t} \rho_0(t,s,X(t)).
%\end{equation*}
for $t$ in $[0,\tau]$ and $s=1,\ldots, B$. Apart from the stratification, this definition of the rate function can be found in~\cite{scheike}.
%pour $s=1,\ldots,B$ puisque pour $s=1$ par exemple on veut estimer la probabilité d'avoir un év. recurrent sachant qu'on n'en a pas encore eu et pour $s=B-1$ on veut estimer la proba. d'avoir un $B$ème év. récurrent sachant qu'on en a eu $B$}. 
%, with cumulative distribution function $F$
 
We consider two semiparametric models for the function $\rho_0$. The first one is an event-specific multiplicative rate model introduced in~\cite{prentice}%, also known as the PWP model
. In this model, the rate function is specified, for $t$ in $[0,\tau]$, by
\begin{align}\label{eq:modelmult}
 \rho_0(t,s,X(t)) = \alpha_0(t,s)\exp\left(X(t)\beta_0(s)\right)
\end{align}
where for each event number $s$, $\beta_0(s)$ is an unknown $p$-dimensional vector of parameters and $\alpha_0$ is an unknown baseline function.

 Following~\cite{scheike}, and \cite{MR2672492}, we also propose to consider its additive counterpart. The rate function in our event-specific additive model is then for $t$ in $[0,\tau]$:
\begin{align}\label{eq:modeladd}
 \rho_0(t,s,X(t))=\left( \alpha_0(t,s)+X(t)\beta_0(s) \right).
\end{align}
%where again, for each event number $s$, $\beta_0(s)$ is an unknown $p$-dimensional vector of parameters.
The models, where $\beta_0$ is constant over the events are refereed to as constant models in what follows.

We consider the problem of estimating the unknown parameter $\beta_0$, in stratified models \eqref{eq:modelmult} and \eqref{eq:modeladd} on the basis of data from $n$ independent and identically distributed random variables. 
Introduce the censoring time  $C$. In a random sample of $n$ subjects, the data consist of $\{N_i(t),T_i,\delta_i,X_i(t), t\leq \tau\}$, $i=1,\ldots,n$ where $N_i(t)=N^*_i(t\wedge C_i)$, $T_i=D_i \wedge C_i$ is the minimum between $D_i$ and $C_i$, $\delta_i=\ind{D_i \leq C_i}$ and $(X_i(t), 0\leq t\leq T_i) $ is the covariates process.
%The data of each subject $i$ ($i=1,\ldots,n$) are denoted by $T_i$, $\delta_i=\ind{D_i \leq C_i}$, $(X_i(t), t\leq T_i) $ and $(N_i(t)=N^*_i(t\wedge C_i), t \leq T_i)$, where $C_i$ is a censoring time and $T_i=D_i \wedge C_i$ is the minimum between $D_i$ and $C_i$. We assume that the censoring times are independent and fulfill Assumption~\ref{ass:independence}.
The next assumption characterizes the dependence mechanism between the censoring time and the other variables.
%We assume that the censoring times fulfill Assumption~\ref{ass:independence}.
\begin{assumption}\label{ass:independence}
For all $s=1,\ldots, B$ and $t$ in $[0,\tau]$, %$i=1,\ldots,n$
\begin{equation*} \E \big( dN^*(t)\mid X(t), D \wedge C \geq  t,N^*(t)=s-1\big)= \E \big( dN^*(t)\mid X(t), D  \geq  t,N^*(t)=s-1\big).\end{equation*}
 \end{assumption}
Note that this assumption is slightly weaker than assuming the independence between $C$ and $(N^*,D, X)$. A similar assumption can be found for instance in~\cite{linwei00}. We also impose the following conditions on the tails of the distribution of $T$ and $N$.
%This assumption is classical in recurrent events context and can be found for instance in~\cite{linwei00}. We also impose the following conditions on the tails of the distribution of $T$ and $N$. and is consequently more general
 \begin{assumption}\label{ass:bounded}
There exists a nonnegative integer $B$ such that
\begin{enumerate}[(i)]
\item $\forall t\in [0,\tau]$, $\mathbb P\big(N(t)\leq B\big)=1$,
\item $\forall t\in [0,\tau]$, $\forall s=1,\ldots,B$, $\mathbb P\big(T\geq t,N(t)=s-1 \mid X(t)\big)>0$.
\end{enumerate} 
\end{assumption}
%Assumption~\ref{ass:bounded} (ii) is used to overcome classical problems in the tail of the distribution of $T$ a classical assumption 
%\texttt{L'Assumption}~\ref{ass:bounded} (ii) \texttt{nous assure que} $s^{(0)}$ \texttt{dans le modèle de Cox est bien strictement positif et que }$\mathbb E[Y^s(t)]>0$ \texttt{dans le modèle de Aalen !!}\\ 
%~\ref{ass:bounded} (i)
Assumption~\ref{ass:bounded} (i) ensures that in models \eqref{eq:modelmult} and \eqref{eq:modeladd}, the total number of observed events is almost surely bounded.
It is standard for inference for recurrent events process, see e.g.~\cite{Dauxois09},~\cite{scheike} or~\cite{BCG}.

Under Assumption~\ref{ass:bounded}, the unknown vector of parameters $\beta_0$ has $p\times B$ unknown coefficients to be estimated. For reasonable sizes of sample $n$, these models are over-parametrized in the sense that, when $\sqrt n \leq p\times B$, the estimators show very poor behaviours (see Section~\ref{sec:simu} for an illustration). On the other hand, simpler forms of models~\eqref{eq:modelmult} and \eqref{eq:modeladd}, in which the unknown parameter does not change with the event,  $\beta_0(s)=\beta_0$, might be too poor to accurately fit the data (see also Section~\ref{sec:simu} and the discussion in~\cite{kelly}). In this paper, we aim at providing estimators realizing a compromise between these two situations.

In the following, we define, for each individual $i$, the event-specific at-risk function $Y_i^s$ and the overall at-risk function $Y_i$ for all $t$ in $[0,\tau]$: 
\begin{align*}
&Y_i^s(t) =\ind{T_i \geq t,N_i(t)=s},\; \;\;Y_i(t)=\sum_{s=1}^BY_i^s(t)=\ind{T_i \geq t}.%{s \in \{1,\ldots,B\}}
\end{align*}

\subsection{Inference in the multiplicative model}\label{par:detailmult}
As in~\cite{prentice}, in the multiplicative event-specific model~\eqref{eq:modelmult}, an estimator $\hat \beta_{ES/mult}$ of the unknown parameter $\beta_0 \in \mathbb R^{p\times B}$ is defined as the maximizer of the partial log-likelihood, or equivalently as
\begin{align}\label{eq:contrastMult}
&\hat \beta_{ES/mult} \in \underset{\beta \in \mathbb R^{p\times B} } \argmin \;L_n^{PL}(\beta)%=  \underset{\beta \in \mathbb R^{p\times B} }\argmin  \left\{ -\frac{1}{n}\sum_{s=1}^B \mathcal L_n(\beta,s)\right\}
\\&=\underset{\beta \in \mathbb R^{p\times B} }\argmin\left[ -\frac{1}{n}\sum_{s=1}^B \sum_{i=1}^n  \int  \left\{X_i(t) \beta(s)-\log\left( \sum_{j=1}^nY_j^s(t)\exp{(X_j(t) \beta(s))} \right)\right\} Y_i^s(t)dN_i(t)   \right].
\nonumber \end{align}
An estimator $\hat \beta_{C/mult}$ in the constant model is defined as
\begin{equation}\label{eq:defConstmult}
\hat \beta_{C/mult} \in  \underset{\beta \in \mathbb R^{p} }\argmin \left[ - \frac{1}{n}\sum_{i=1}^n  \int  \left\{X_i(t)\beta -\log\left( \sum_{j=1}^nY_j(t)\exp{(X_j(t) \beta)} \right)\right\} Y_i(t)dN_i(t)   \right].
\end{equation}

\subsection{Inference in the additive model}\label{par:detailadd}
As noticed in~\cite{martinussenscheike09a,martinussenscheike09b} or~\cite{GGaalen}, in the usual additive hazards model, the estimator $\hat \beta_{ES/add}$ of the unknown parameter $\beta_0 \in \mathbb R^{p\times B}$ can be written as the minimizer of a (partial) least-squares criterion:
%\texttt{ vérifier les matrices après avoir écrit la preuve de la proposition1 !!!!!!!!!!!!!!!!!}
\begin{align}\label{eq:contrastAdd}
&\hat \beta_{ES/add} \in  \underset{\beta \in \mathbb R^{p\times B} } \argmin \;L_n^{PLS}(\beta)= \underset{\beta \in \mathbb R^{p\times B} }\argmin \sum_{s=1}^B \left\{ \beta(s)^\top \H_n(s)\beta(s) -2 \bs h_n(s) \beta(s)  \right\},
\end{align}
where for all $s \in \{1,\ldots,B\}$, $\H_n(s)$ are $p \times p$ symetrical positive semidefinite
matrices equal to
\begin{equation*}
\frac 1n \sum_{i=1}^n \int Y_i^s(t) \Big(X_i(t) -\bar X^s(t)\Big)^{\otimes 2}dt,
\end{equation*}
and where ${\bs h}_n(s)$ are $p$-dimensional vectors equal to
\begin{equation*}
\frac 1n \sum_{i=1}^n  \int\ind{N_i(t)=s}   \Big(X_i(t)- \bar X^s(t)\Big)dN_i(t),
\end{equation*}
with $\bar X^s(t)=\sum_{i=1}^n X_i(t)Y_i^s(t)/\sum_{i=1}^nY_i^s(t)$.
We show in the Appendix why this criterion is a relevant strategy in the additive event-specific model.

On the other hand, an estimator $\hat \beta_{C/add}$ in the constant model is defined as
\begin{align}\label{eq:defConstadd}
&\hat \beta_{C/add} \in  \underset{\beta \in \mathbb R^{p} }\argmin \left(  \beta^\top  \H_n \beta -2 \bs h_n \beta \right), \text{ with } \H_n= \sum_{s=1}^B \H_n(s) \text{ and }\bs h_n=\sum_{s=1}^B  \bs h_n(s).
\end{align}
%with $\H_n= \sum_{s=1}^B \H_n(s)$ and $\bs h_n=\sum_{s=1}^B  \bs h_n(s)$.

\subsection{A total-variation penalty}\label{par:TVpen}
To overcome the possible over-parametrization of models~\eqref{eq:modelmult} and~\eqref{eq:modeladd}, we propose to define penalized versions of criteria~\eqref{eq:contrastMult} and~\eqref{eq:contrastAdd}. %Typical behaviours of the three estimators, are illustrated on a simulation study in Section~\ref{sec:simu}.
For all $\beta=(\beta(s), s=1,\ldots,B)$ with $\beta(s)=(\beta^1(s),\ldots,\beta^p(s))$, define for all $j=1,\ldots,p $
\begin{align}\label{eq:betaj}
\beta^j=(\beta^j(1),\ldots,\beta^j(B))  \text{ and } 
\textsc{tv}(\beta^j)=\sum_{s=2}^B |\beta^j(s)-\beta^j(s-1)|= \sum_{s=2}^B|\Delta \beta^j(s)|.%\nonumber
\end{align}

We now consider the minimizers of the partial log-likelihood (respectively the partial least-squares) penalized with a covariate specific total variation. Define  the penalized estimators in models~\eqref{eq:modelmult} and ~\eqref{eq:modeladd} as:
\begin{align}\label{eq:defTVmult}
\hat \beta_{\textsc{tv}/mult} &\in\underset{\beta \in \mathbb R^{p\times B} }{\argmin} \; \left\{L_n^{PL}(\beta)+\frac{\lambda_n}{n} \sum_{j=1}^p \textsc{tv}(\beta^j)\right\} \text{ and}\end{align}
\begin{equation}\label{eq:defTVadd}
\hat \beta_{\textsc{tv}/add} \in \underset{\beta \in \mathbb R^{p\times B} }{\argmin} \; \left\{L_n^{PLS}(\beta)+\frac{\lambda_n}{n}\sum_{j=1}^p \TV(\beta^j)\right\}.
\end{equation}
%where
%\begin{align*}
%M_n^{mult}(\beta) & = L^{PL}(\beta)+\frac{\lambda_n}{n} \sum_{j=1}^p \sum_{s=2}^B | \beta^j(s)-\beta^j(s-1)|,\\
%M_n^{add}(\beta) & =  L^{PLS}(\beta)+\frac{\lambda_n}{n}\sum_{j=1}^p \sum_{s=2}^B | \beta^j(s)-\beta^j(s-1)|.
%\end{align*}
These penalized algorithms can be rewritten as lasso algorithms (the details are given in Supplementary Material).
%replicates $(N_i,T_i,X_i,\delta_i)_{i=1,\dots,n}$ of $(N,T,X,\delta)$.
%In the following we will respectively denote by $H$ and $F_X$ the cumulative distribution functions of $T$ and $X$.  
% \begin{assum}\label{ass:indep}
%Assume that:
%\begin{enumerate}[(i)]
%\item $\P \big[dN^*(C)\neq 0\big]=0$,
%\item $\P[D=C]=0,$
%\item $C$ is independent of $D$ and of the process $(N(t))_{t\geq 0}$,
%\item $\P[C\leq t\midN^*,X,D]=\P[C\leq t\midN^*,D]=G(t) \text{ for } t\in [0,\tau_H]$.
%\end{enumerate}
%\end{assum}
\section{Asymptotic results}\label{sec:asymp}
We successively provide the asymptotic results for the estimators $\hat \beta_{\textsc{tv}/add}$ in the additive model and $\hat \beta_{\textsc{tv}/mult}$ in the multiplicative model. In both models, the following condition is mandatory.\begin{assumption}\label{ass:covbounded}
The covariates process $X(\cdot)$ is of bounded variation on $[0,\tau]$.
\end{assumption}
Define for all $s=1,\ldots, B$ the centered process $M^s(t)=  N(t)-\E \big(N(t)\mid X(t), D \wedge C \geq  t,N(t)=s-1\big)$ and the $p\times p$ matrix 
%Define for all $s=1,\ldots, B$, the following $p\times p$ matrix
\[\H(s) :=\int \mathbb E[Y^s(t)X(t)^\top X(t)]dt-\int\frac{(\mathbb E[Y^s(t)X(t)])^{\otimes 2}}{\mathbb E[Y^s(t)]}dt,\]
which from Assumption~\ref{ass:bounded} (ii) is well defined. 
\begin{theorem}\label{theo:CVadd}
%Assume that, for each $s=1,\ldots,B$, $\H_n(s)$ converges to a non-singular matrix $\H(s)$ and that Asumptions~\ref{ass:bounded}, \ref{ass:independence} and~\ref{ass:covbounded}
%are fulfilled.
Assume that, for each $s=1,\ldots,B$, $\H(s)$ is non-singular and that Asumptions~\ref{ass:independence},~\ref{ass:bounded} and~\ref{ass:covbounded}
are fulfilled.

1. If $\lambda_n/n \to 0$ as $n\to \infty$ then $\hat \beta_{\textsc{tv}/add}$ converges to $\beta_0$ in probability.

2.  If $\lambda_n/\sqrt n \to \lambda_0\geq 0$ as $n\to \infty$ then $\sqrt n (\hat \beta_{\textsc{tv}/add}-\beta_0)$ converges in distribution to 
\begin{align*}
& \underset{u \in \mathbb R^p}{\argmin} \;\Lambda_{add}(u) = \underset{u \in \mathbb R^p}{\argmin} \; \Big[\sum_{s=1}^B \left\{u(s)^\top \H(s) u(s)-2u(s)^\top\xi_{add}(s)\right\}
\\& 
+\lambda_0\sum_{j=1}^p \sum_{s=2}^B \left\{| \Delta u^j(s)|\ind{\Delta \beta^j(s)=0}\right.
%\\&
 %\qquad\quad 
 +\left.\text{sgn}(\Delta \beta^j(s))(\Delta u^j(s))\ind{\Delta \beta^j(s)\neq0}\right\} \Big],
\end{align*}
and for each $s$, $\xi_{add}(s)$ is a centered p-dimensional gaussian vector with covariance matrix equal to
\[\mathbb E\left[\Big(\int_0^{\tau}(X(t)-\mathbb E[Y^s(t)X(t)]/\mathbb E[Y^s(t)])\ind{N(t)=s}dM^s(t)\Big)^{\otimes 2}\right].\]

\end{theorem}

%The asymptotic results for the Cox model can be proved under classical conditions on the model, such as conditions VII.2.1 and VII.2.2 pages 497-498 and (7.2.1) page 481 in~\cite{ABGK}. However, our Assumptions~\ref{ass:bounded} and~\ref{ass:covbounded} provide sufficient conditions for these to hold (see for instance example VII.2.7 page 502 in~\cite{ABGK}).

%Therefore, d
Define for all $s=1,\ldots, B$ and for all $t\in [0,\tau]$,
\[s^{(l)}(s,t,\beta)=\mathbb E[Y^s(t)X(t)^{\otimes l}\exp(X(t) \beta(s))], l=0,1,2.\]
Introduce $\mathbf e(s,t,\beta)=s^{(1)}(s,t,\beta)/s^{(0)}(s,t,\beta)$, $\mathbf v(s,t,\beta)=s^{(2)}(s,t,\beta)/s^{(0)}(s,t,\beta)-\mathbf e(s,t,\beta)^{\otimes 2}$ and $\mathbf \Sigma(s,\beta)=\int \mathbf v(s,t,\beta)\mathbb E[Y^s(t)dN(t)]$. For any $s=1,\ldots, B$ and for any $t\in [0,\tau]$, the three functions $s^{(l)}(s,t,\beta_0)$ are bounded from Assumption~\ref{ass:covbounded} and $\mathbf e(s,t,\beta), \mathbf v(s,t,\beta)$ and $\mathbf \Sigma(s,\beta)$ are finite from Assumptions~\ref{ass:bounded} and~\ref{ass:covbounded}. 

\begin{theorem}\label{theo:CVmult}
Assume that for each $s=1,\ldots,B$, $\mathbf \Sigma(s,\beta_0)$ is non-singular and that Assumptions~\ref{ass:independence},~\ref{ass:bounded} and~\ref{ass:covbounded} are fulfilled. 

1. If $\lambda_n/n \to 0$ as $n\to \infty$ then $\hat \beta_{\textsc{tv}/mult}$ converges to $\beta_0$ in probability.

2.  If $\lambda_n/\sqrt n \to \lambda_0\geq 0$ as $n\to \infty$ then $\sqrt n (\hat \beta_{\textsc{tv}/mult}-\beta_0)$ converges in distribution to
% \begin{align*}
% V(u)& =\sum_{s=1}^B \left\{\frac 12 u(s)^\top \int \left(\frac{s^{(2)}(s,t,\beta_0)}{s^{(0)}(s,t,\beta_0)}-\Bigg(\frac{s^{(1)}(s,t,\beta_0)}{s^{(0)}(s,t,\beta_0)}\Bigg)^{\otimes 2}\right)\,\mathbb E[Y^s(t)dN(t)]u(s)\right\}\\
% & \quad +\sum_{s=1}^B u(s)^\top\Sigma(s)+\lambda_0\sum_{j=1}^p \sum_{s=2}^B \left\{| \Delta u^j(s)|\mathds 1_{\beta_0^j(s)=\beta_0^j(s-1)}\right.\\
% & \quad+\left.\text{sgn}(\beta_0^j(s)-\beta_0^j(s-1))(\Delta u^j(s))\mathds 1_{\beta_0^j(s)\neq \beta_0^j(s-1)}\right\},
% \end{align*}
\begin{align*}
& \underset{u \in \mathbb R^p}{\argmin} \;\Lambda_{mult}(u) = \underset{u \in \mathbb R^p}{\argmin} \; \Big[\sum_{s=1}^B \left\{\frac 12 u(s)^\top \mathbf \Sigma(s,t,\beta_0)u(s)+u(s)^\top\mathbf\xi_{mult}(s)\right\}\\
& \quad +\lambda_0\sum_{j=1}^p \sum_{s=2}^B \left\{| \Delta u^j(s)|\ind{\Delta\beta_0^j(s)=0}+\text{sgn}(\Delta \beta_0^j(s))(\Delta u^j(s))\ind{\Delta\beta_0^j(s)\neq0}\right\}\Big],
\end{align*}
and for each $s$, $\xi_{mult}(s)$ is a centered p-dimensional gaussian vector with covariance matrix equal to
\[\mathbb E \left[\left(\int_0^\tau \left(X(t)-\mathbf e(s,t,\beta_0)\right)Y^s(t)dM^s(t)\right)^{\otimes2}\right] .\]
%\end{enumerate}

\end{theorem}

Theorems~\ref{theo:CVadd} and~\ref{theo:CVmult} prove the consistency and asymptotic normality of our estimators~(\ref{eq:defTVmult}) and~(\ref{eq:defTVadd}). This assures that they behave better than the constant estimators when $\beta_0$ is non constant. In addition, the considered penalty will induce sparsity for each covariate $j= 1,\ldots,p$ in the successive differences $\Delta \beta^j(s)$, $s=1,\ldots,B$.  As a consequence, the effects of a covariate on two consecutive events will often be equal.
We show, in the following simulation study, that this induced sparsity ameliorates the behaviour of our estimators compared to the unconstrained ones (defined in Equations~\eqref{eq:contrastMult} and~\eqref{eq:contrastAdd}).

%\texttt{ a changer }The rationale behind the choice of a total-variation penalty is to make use of the chronological nature of the data and automatically estimate a simple model, where 
\section{Simulation studies}\label{sec:simu}
We compare the performances of the penalized estimators \eqref{eq:defTVmult} and~\eqref{eq:defTVadd}, the constant ones~\eqref{eq:defConstmult} and~\eqref{eq:defConstadd} , and the unconstrained ones~\eqref{eq:contrastMult} and~\eqref{eq:contrastAdd}. To mimic the bladder tumour cancer dataset studied in Section~\ref{sec:real}, we set $p=4$ and consider $B=5$ recurrent events for the estimation. In the multiplicative and additive models, the sample size $n$ varies from $n =50= 2\!\cdot\!5\; pB$ to $n=1000 \simeq (pB)^{2\cdot3}$.

We draw the $p=4$ covariates from uniform distributions and set the parameters values at $\beta_0^1=(0,0,b_1,b_1,0,\ldots,0)$, $\beta_0^2=(b_2,\ldots,b_2)$,  $\beta_0^3=b_3(1,2,3,\ldots)$ and $\beta_0^4=(0,\ldots,0)$. We generate recurrent event times from the multiplicative~\eqref{eq:modelmult} and additive~\eqref{eq:modeladd} models with baseline defined through the Weibull distribution with shape parameter $a_{\mathcal W}$ and scale parameter $1$. The death and censoring times are generated from exponential distributions with parameters $a_D$ and $a_C$ respectively. 
%More precisely, we generate independent inter-events times from event specific multiplicative and additive semiparametric models for the intensity. We stop the simulation at the minimum between death and censoring times. This implies the event-specific models for the rate function defined in~\eqref{eq:modelmult} and~\eqref{eq:modeladd}. 
We set the value of parameter $a_{\mathcal W}$ at $2\!\cdot\!5$. Finally, the values of $a_D$ and $a_C$ are empirically determined to obtain $\textsc{p}_{\text{obs}}=28-29\%$ and $14-15\%$ of individuals experiencing the fifth event.

%The constant and event specific estimators are calculated via the functions \texttt{coxph} (in the R package \texttt{survival}, see~\cite{survival}) and \texttt{ahaz} (in the R package \texttt{ahaz}, see~\cite{ahaz}). We develop ad hoc functions to calculate the penalized estimators, they are based on the functions \texttt{coxnet} (in the R package \texttt{glmnet}, see~\cite{glmnet}) and \texttt{ahazpen} (in the R package \texttt{ahaz})% (see Supplementary Material for details)
%. In particular, we choose the regularization parameters via cross-validation, as proposed in the latter packages.

To evaluate the performances of the different estimators, we conduct a Monte Carlo study with $M=200$ experiences. The estimation accuracy is investigated for each method via a mean squared rescaled error defined as
\begin{equation}\label{eq:mse}
\textsc{mse}=\frac{1}{M}\sum_{m=1}^M \frac{\norm{\hat \beta_m-\beta_0}^2}{\norm{\beta_0}^2},
\end{equation}where $\hat \beta_m$ is the estimation in the sample $m$. We furthermore study the detection power of non-constant (respectively constant) covariate effects by computing mean false positive (\textsc{fp}) rates and mean false negative (\textsc{fp}) rates for each method. They are defined, for an estimation $\hat \beta_m$, as
\begin{equation}\label{eq:FP}
\textsc{fp}(\hat \beta_m)=\card\left(j \in \{1,\ldots,p\} \text{ s.t. } \TV(\hat \beta^j) \neq 0  \text{ and }  \TV(\beta_0^j) = 0\right)\end{equation}and
\begin{equation}\label{eq:FN}
\textsc{fn}(\hat \beta_m)=\card\left(j \in \{1,\ldots,p\} \text{ s.t. } \TV(\hat \beta^j) = 0   \text{ and }\TV(\beta_0^j) \neq 0\right),
\end{equation}
where $\TV$ is defined in~\eqref{eq:betaj}.

\begin{table}
\def~{\hphantom{0}}
\tbl{Simulation results in the multiplicative model for $\textsc{p}_{\text{obs}}=28\%$}{%
\begin{tabular}{c||ccc|ccc|ccc|ccc}
n&\multicolumn{3}{c|}{Unconstrained}&\multicolumn{3}{c|}{Constant}&\multicolumn{3}{c|}{\textsc{tv}}&\multicolumn{3}{c}{two-steps \textsc{tv}}\\
&\textsc{mse}&\textsc{fp}&\textsc{fn}&\textsc{mse}&\textsc{fp}&\textsc{fn}&\textsc{mse}&\textsc{fp}&\textsc{fn}&\textsc{mse}&\textsc{fp}&\textsc{fn}\\\hline\hline
50 &    0$\cdot$100    &       2       &     0   &  0$\cdot$412  &   0&              2       &   0$\cdot$054 &     1$\cdot$44       &    0$\cdot$03   &     0$\cdot$044    &     0$\cdot$82        &     0$\cdot$02\\     
  100&    0$\cdot$030     &       2       &     0  &  0$\cdot$415&   0&              2     &   0$\cdot$025  &    1$\cdot$54        &   0  &    0$\cdot$019    &      0$\cdot$76      &       0\\
500    &0$\cdot$006      &       2       &     0  & 0$\cdot$413  &   0&              2           & 0$\cdot$008      & 1$\cdot$76          & 0   & 0$\cdot$006     &       0$\cdot$30    &      0\\
1000  & 0$\cdot$005     &       2       &     0   & 0$\cdot$415&   0&              2 & 0$\cdot$006     & 1$\cdot$81          & 0     &  0$\cdot$006               &           0$\cdot$05           & 0\\
 \end{tabular}
 }
\label{table:mult28}
\begin{tabnote}
\textsc{mse}: mean squared error, \textsc{fp}: false positives,  \textsc{fn}: false negatives.
\end{tabnote}
\end{table}
%
%     mean(err_PWP) mean(err_Const) mean(err_\textsc{tv}_min) mean(err_\textsc{tv}_adap) mean(obs) mean(cens) mean(FP_PWP) mean(FN_PWP) mean(FP_\textsc{tv}_min) mean(FN_\textsc{tv}_min) mean(FP_Const)
%50     0.100050010       0.4122113       0.05439090       0.046282833   0.28100   0.254000            2            0           1.445           0.025              0
%100    0.030416200       0.4153868       0.02519778       0.019270662   0.28565   0.251000            2            0           1.540           0.000              0
%500    0.006316011       0.4134805       0.00767808       0.006458248   0.28677   0.250750            2            0           1.755           0.000              0
%1000   0.004474525       0.4147892       0.00555721       0.006197266   0.28373   0.250365            2            0           1.805           0.000              0
%     mean(FN_Const) mean(FP_\textsc{tv}_adap) mean(FN_\textsc{tv}_adap)
%50                2             0.73             0.02
%100               2             0.76             0.00
%500               2             0.30             0.00
%1000              2             0.05             0.00

%[1] "Simu pour a= 1.2 c = 0.4 a0 = 2.5 model mult"

\begin{table}
\def~{\hphantom{0}}
\tbl{Simulation results in the multiplicative model for $\textsc{p}_{\text{obs}}=14\%$}{%
\begin{tabular}{c||ccc|ccc|ccc|ccc}
n&\multicolumn{3}{c|}{Unconstrained}&\multicolumn{3}{c|}{Constant}&\multicolumn{3}{c|}{\textsc{tv}}&\multicolumn{3}{c}{two-steps \textsc{tv}}\\
&\textsc{mse}&\textsc{fp}&\textsc{fn}&\textsc{mse}&\textsc{fp}&\textsc{fn}&\textsc{mse}&\textsc{fp}&\textsc{fn}&\textsc{mse}&\textsc{fp}&\textsc{fn}\\\hline\hline
50&\textsc{NA}&\textsc{NA}&\textsc{NA}&0$\cdot$440 &0&2&     0$\cdot$161   &  1$\cdot$37 &          0$\cdot$185      &      0$\cdot$137& 0$\cdot$82  &           0$\cdot$19 \\
100 &   0$\cdot$566&  2 &           0        &0$\cdot$434&0&2&     0$\cdot$053&1$\cdot$55      &     0$\cdot$005  &         0$\cdot$042&  0$\cdot$88       &         0\\   
500    &0$\cdot$014&  2 &           0        &       0$\cdot$433&0&2&  0$\cdot$016       &  1$\cdot$84   &        0   &      0$\cdot$012&           1$\cdot$06    &           0 \\
1000   &0$\cdot$009     &  2 &           0        & 0$\cdot$433    &0&2&  0$\cdot$011     &    1$\cdot$89        & 0       &  0$\cdot$010&   0$\cdot$68    &            0 \\  
 \end{tabular}
 }\label{table:mult14}
\begin{tabnote}
\textsc{mse}: mean squared error, \textsc{fp}: false positives,  \textsc{fn}: false negatives,  \textsc{na}: non applicable .
\end{tabnote}
\end{table}
%
%[1] "Simu pour a= 1.2 c = 0.4 a0 = 2.5 model mult"
%     mean(err_PWP) mean(err_Const) mean(err_\textsc{tv}_min) mean(err_\textsc{tv}_adap) mean(obs) mean(cens) mean(FP_PWP) mean(FN_PWP) mean(FP_\textsc{tv}_min) mean(FN_\textsc{tv}_min) mean(FP_Const)
%100    0.566179461       0.4336312       0.05337511        0.04150011   0.14240    0.24635            2            0            1.55           0.005              0
%500    0.013476021       0.4328611       0.01571096        0.01219624   0.14316    0.24968            2            0            1.84           0.000              0
%1000   0.009447454       0.4329942       0.01136878        0.01016262   0.14238    0.24967            2            0            1.89           0.000              0
%     mean(FN_Const) mean(FP_\textsc{tv}_adap) mean(FN_\textsc{tv}_adap)
%100               2            0.875                0
%500               2            1.055                0
%1000              2            0.675                0
%
%[1] "Simu pour a= 1.2 c = 0.4 a0 = 2.5 model mult"
  % mean(err_PWP) mean(err_Const) mean(err_\textsc{tv}_min) mean(err_\textsc{tv}_adap) mean(obs) mean(cens) mean(FP_PWP) mean(FN_PWP) mean(FP_\textsc{tv}_min) mean(FN_\textsc{tv}_min) mean(FP_Const)
%50             0       0.4404453        0.1608706         0.1639882    0.1448     0.2537            0            0            1.37           0.185              0
   %mean(FN_Const) mean(FP_\textsc{tv}_adap) mean(FN_\textsc{tv}_adap)
%50              2            0.855             0.19

%[1] "Simu pour a= 3.6 c = 1.2 a0 = 2.5 model add"

\begin{table}
\def~{\hphantom{0}}
\tbl{Simulation results in the additive model for $\textsc{p}_{\text{obs}}\simeq28\%$}{%
\begin{tabular}{c||ccc|ccc|ccc|ccc}
n&\multicolumn{3}{c|}{Unconstrained}&\multicolumn{3}{c|}{Constant}&\multicolumn{3}{c|}{\textsc{tv}}&\multicolumn{3}{c}{two-steps \textsc{tv}}\\
&\textsc{mse}&\textsc{fp}&\textsc{fn}&\textsc{mse}&\textsc{fp}&\textsc{fn}&\textsc{mse}&\textsc{fp}&\textsc{fn}&\textsc{mse}&\textsc{fp}&\textsc{fn}\\\hline\hline
50   &    4$\cdot$986    & 2         &   0   &   0$\cdot$416   & 0         &  2   &    0$\cdot$467     &0$\cdot$98     &      0$\cdot$58   & 1$\cdot$142&             0$\cdot$65      &      0$\cdot$81  \\
100&      0$\cdot$935   & 2         &   0   &    0$\cdot$351   & 0         &  2   &         0$\cdot$254      &      1$\cdot$38    &       0$\cdot$21 &  0$\cdot$353& 0$\cdot$86     &       0$\cdot$48 \\
500   &   0$\cdot$135     & 2         &   0   & 0$\cdot$309   & 0         &  2   &         0$\cdot$079     &     1$\cdot$91         &  0$\cdot$01 &   0$\cdot$094&1$\cdot$44      &      0$\cdot$08   \\
1000   &  0$\cdot$071   & 2         &   0   &   0$\cdot$299     & 0         &  2   &       0$\cdot$049  &   1$\cdot$98      &     0 &     0$\cdot$05& 1$\cdot$64      &      0 \\
 \end{tabular}}\label{table:add28}
\begin{tabnote}
\textsc{mse}: mean squared error, \textsc{fp}: false positives,  \textsc{fn}: false negatives
\end{tabnote}
\end{table}
%
%     mean(err_PWP) mean(err_Const) mean(err_\textsc{tv}_min) mean(err_\textsc{tv}_adap) mean(obs) mean(cens) mean(FP_PWP) mean(FN_PWP) mean(FP_\textsc{tv}_min) mean(FN_\textsc{tv}_min) mean(FP_Const)
%50       4.9856779       0.4164235       0.46669071        0.68758869   0.28570    0.25070            2            0           0.980           0.575              0
%100      0.9353329       0.3506671       0.25394102        0.35338775   0.28655    0.24870            2            0           1.380           0.205              0
%500      0.1353786       0.3089726       0.07927459        0.09438643   0.28866    0.24896            2            0           1.905           0.010              0
%1000     0.0711683       0.2992331       0.04931978        0.05419548   0.28968    0.24851            2            0           1.975           0.000              0
%     mean(FN_Const) mean(FP_\textsc{tv}_adap) mean(FN_\textsc{tv}_adap)
%50                2             0.61            0.780
%100               2             0.86            0.475
%500               2             1.44            0.080
%1000              2             1.64            0.000
%[1] "Simu pour a= 6.6 c = 2.2 a0 = 2.5 model add"

\begin{table}
\def~{\hphantom{0}}
\tbl{Simulation results in the additive model for $\textsc{p}_{\text{obs}}\simeq14\%$}{%
\begin{tabular}{c||ccc|ccc|ccc|ccc}
n&\multicolumn{3}{c|}{Unconstrained}&\multicolumn{3}{c|}{Constant}&\multicolumn{3}{c|}{\textsc{tv}}&\multicolumn{3}{c}{two-steps \textsc{tv}}\\
&\textsc{mse}&\textsc{fp}&\textsc{fn}&\textsc{mse}&\textsc{fp}&\textsc{fn}&\textsc{mse}&\textsc{fp}&\textsc{fn}&\textsc{mse}&\textsc{fp}&\textsc{fn}\\\hline\hline
50    &      \textsc{NA} &\textsc{NA}&\textsc{NA}&      0$\cdot$505   &0&2&   0$\cdot$781      &0$\cdot$95     &      0$\cdot$81  &  2$\cdot$368 &   0$\cdot$86       &     0$\cdot$97\\ 
100     & 4$\cdot$114   &2&0&    0$\cdot$393  &0&2&      0$\cdot$707     &1$\cdot$450   &        0$\cdot$27 &    0$\cdot$84 & 1$\cdot$11      &      0$\cdot$52\\ 
500     & 0$\cdot$339   &2&0&    0$\cdot$330   &0&2&     0$\cdot$154    &1$\cdot$975      &     0$\cdot$01   &     0$\cdot$19 &1$\cdot$67   &         0$\cdot$06\\
1000    & 0$\cdot$171  &2&0&     0$\cdot$320   &0&2&     0$\cdot$097      &  1$\cdot$995   &        0&   0$\cdot$12  &    1$\cdot$80 &           0$\cdot$02\\
 \end{tabular}}\label{table:add14}
\begin{tabnote}
\textsc{mse}: mean squared error, \textsc{fp}: false positives,  \textsc{fn}: false negatives,  \textsc{na}: non applicable .
\end{tabnote}
\end{table}
%     mean(err_PWP) mean(err_Const) mean(err_\textsc{tv}_min) mean(err_\textsc{tv}_adap) mean(obs) mean(cens) mean(FP_PWP) mean(FN_PWP) mean(FP_\textsc{tv}_min) mean(FN_\textsc{tv}_min) mean(FP_Const)
%50              NA       0.5051657       0.78137314         1.3682388  0.142700    0.24660           NA           NA           0.945           0.805              0
%100      4.1138576       0.3934154       0.70688247         0.8419203  0.145800    0.24880            2            0           1.450           0.270              0
%500      0.3388546       0.3296991       0.15351121         0.1920786  0.138560    0.25192            2            0           1.975           0.005              0
%1000     0.1712242       0.3196165       0.09731826         0.1210118  0.140565    0.25215            2            0           1.995           0.000              0
%     mean(FN_Const) mean(FP_\textsc{tv}_adap) mean(FN_\textsc{tv}_adap)
%50                2            0.775            0.950
%100               2            1.110            0.515
%500               2            1.670            0.060
%1000              2            1.795            0.020

As expected, the constant model is biased and behave poorly for our choice of a non-constant $\beta_0$. The comparison between the unconstrained and penalized estimators is in favour of our estimator in all four cases as long as $n$ is smaller than $p^2$. When the percentage of individuals experiencing the fifth event drops, non-constant estimators are slightly less accurate. Algorithms are not able to compute all $M=200$ unconstrained estimators for $n=50$. % resulting with non applicable computation of mean squared error, false positive and negative rates. 
For $p=4$, $B=5$, $n=100$ and $\textsc{p}_{\text{obs}}=14\%$ (which are values close to those encountered in the bladder tumour cancer dataset studied in the next section) our penalized estimators are respectively $5\!\cdot\!8$, in the additive model, and $10\!\cdot\!6$, in the multiplicative model, times better than the unconstrained ones in terms of estimation error.

Surprisingly the number of false positives detected by our penalized estimators increases when the sample size increases. 
%This is, in our opinion, due to the choice of the regularization parameter performed by packages \texttt{glmnet} and \texttt{ahazpen}. In both packages, the measures used in the cross-validation (respectively the partial log-likelihood and partial leas-squares) are tuned to minimize the estimation error and not the selection accuracy. 
A possible solution to ameliorate the latter is to apply the reweighed lasso, or two-steps lasso, as proposed in~\cite{reweightedlasso}  (details are given in Supplementary Material). We compute the mean squared error, false positive and negative rates of the resulting estimator. It shows better false positive rates than the first step penalized estimator, greater false negative rates and comparable mean squared errors.

We repeat the simulation study for $a_{\mathcal W}=2\!\cdot\!5$ and then for a Gompertz baseline with shape parameter $a_{\mathcal G}=0\!\cdot\!5$ (and $a_{\mathcal G}=0\!\cdot\!5$) and scale parameter $1$.
 The results are reported in Supplementary Material. Conclusions are similar.
 
\section{Bladder tumour data analysis}\label{sec:real}
In this section we illustrate the behaviour of our estimators on the bladder tumour cancer data of~\cite{byar}. These data were obtained from a clinical trial conducted by the Veterans Administration Co-operative Urological Group. One hundred and sixteen patients were randomised to one of three treatments: placebo, pyridoxine or thiotepa. For each patient, the time of recurrence tumours were recorded until the death or censoring times. The number of recurrences ranges from $0$ to $10$. %Two patients who were censored before the beginning of the study are removed.
%%Two outliers were removed from the study. 
On the $n=116$ %remaining
patients,  since $13\!\cdot\!79\%$ experienced at least five tumour recurrences and only $6\!\cdot\!9\%$ patients experienced six tumour recurrences or more, we set the parameter $B$ to $5$. In addition to the two treatment variables, pyridoxine and thiotepa, two supplementary covariates were recorded for each patient: the number of initial tumours and the size of the largest initial tumour. 
%On the $106$ remaining patients,  $53.45\%$ experienced at least one tumour recurrence, $33.62\%$ experienced at least two tumour recurrences, $24.14\%$ experienced at least three tumour recurrences, $17.24\%$ experienced at least four tumour recurrences and $13.79\%$ experienced at least five tumour recurrences. Since only $6.9\%$ patients experienced six tumour recurrences or more, these recurrences are not considered in the estimation procedure. 
%In addition, two supplementary covariates were recorded for each patient: the number of initial tumours and the size of the largest initial tumour. For interpretation purpose, the treatment variable is also coded as two new binary variables, pyridoxine and thiotepa, making placebo the reference. In our setting, with the previous notations, $p=4$, $B=5$, $n=116$.

\begin{figure}[htbp]
% The arguments in the next line are {height}{optional width}{used only by OUP for typesetting}[filename, in directory art]
\includegraphics[height=3.8in,width=6in]{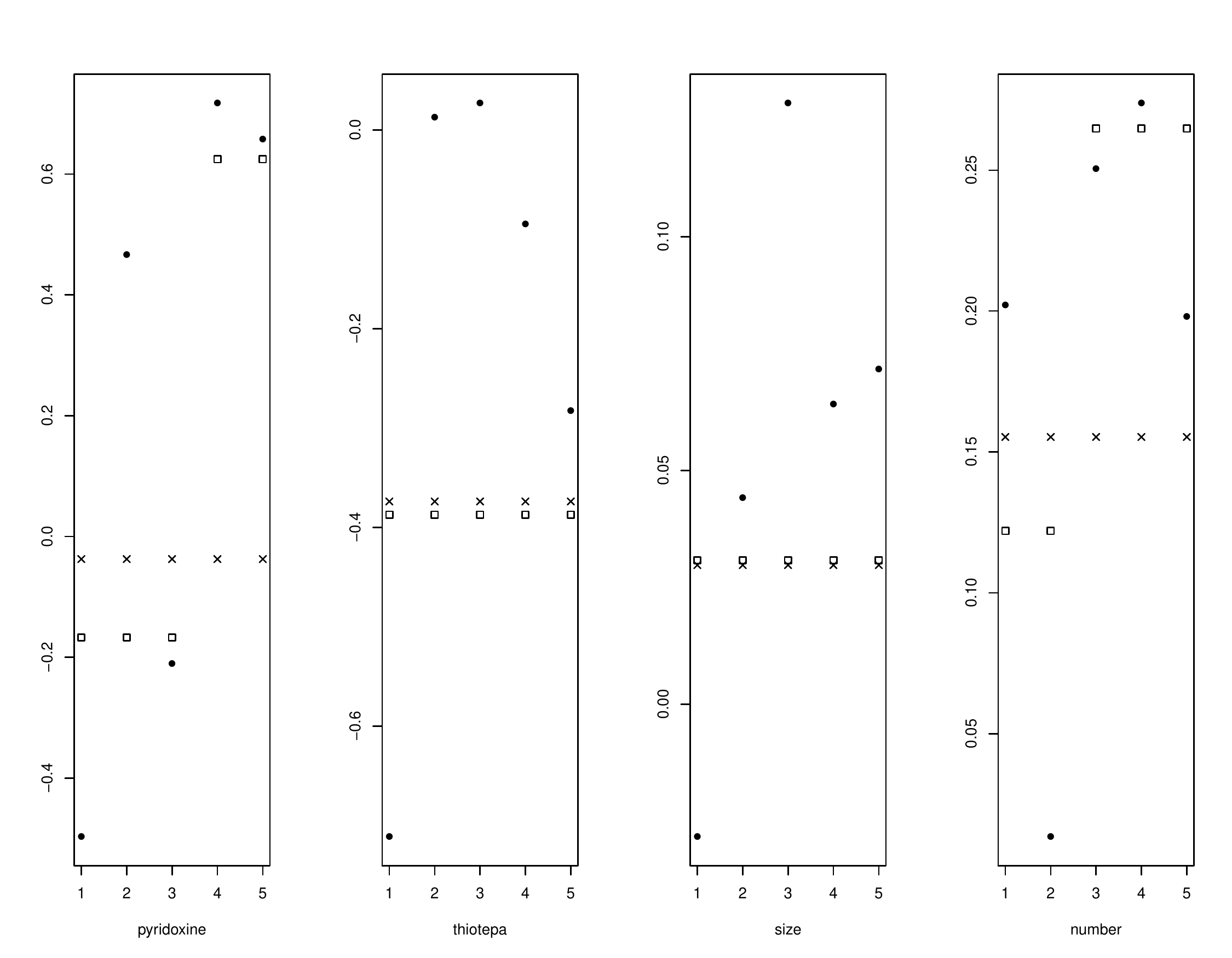} 
%\figurebox{20pc}{25pc}{AdaptMult.pdf]}[AdaptMult.pdf] %AdaptMult.pdf
% note that files may not be rotated
\caption{Estimates for the bladder data in the multiplicative model. The crosses represent the constant estimator, the filled circles the unconstrained estimator and the squares the reweighed lasso estimator.}
\label{fig1}
\end{figure}

Figure~\ref{fig1} displays the estimations obtained from the constant, unconstrained and total variation estimators in the multiplicative model. In order to enforce the variables selection performance of the total variation estimator, the coefficients were estimated using the reweighed lasso. The unconstrained estimator shows very strong variations and is difficult to interpret as such. On the other hand, the constant estimator gives valuable information on the impact of each covariate, but in turn cannot detect a change in variation. Our total-variation estimator reaches compromise: it is not constant but easily interpretable. %between the constant and pwp estimators: it gives a more accurate estimation and is easily interpretable. 
%A nice feature of these estimators is also the ability to detect a possible variation in the effect of the covariates with respect to the number of recurrent events experienced by a subject. 

For instance, a remarkable aspect of the pyrodixine treatment can be highlighted from the total variation estimation: this treatment produces a protective effect for the first three %established, noticed, highlighted
%second and third 
tumour recurrences but
 %as soon as a patient has already experienced three tumours recurrences,
  %the pyrodixine effect is reversed and
   the odds of further recurrences are increased by this treatment. In the same way, an increase in the effect of the initial number of tumours on recurrences is observed from the third recurrence. On the opposite, the effects of the thiotepa treatment or the size of the largest tumour are shown to be constant in the total variation model, the parameter estimates having values similar to the ones obtained in the constant model.

Our conclusions on the treatments effects are in agreement with previous studies on bladder tumours recurrences. For instance, no difference in the rate or time to tumour recurrence was found from patients using pyrodixine with patients using placebo  in~\cite{JOncol} and~\cite{Uro}. Moreover,~\cite{huang03} and~\cite{sunetal06} have respectively studied gap time recurrences in the multiplicative and additive models. The results obtained from the former showed a small protective effect of this treatment while the latter concluded that gap times did not seem related to pyridoxine. 
These examples illustrate the nice features of our total-variation estimator: it provides sharper results, giving relevant informations on covariates effect with respect to the number of recurrent events experienced by a subject and it provides the ability to detect a change of variation. Further details are provided in Supplementary Material.

\section{Discussion}\label{sec:discussion}

In this paper, the Aalen and Cox models were studied to model the effect of covariates on the rate function. However, such models are not essential in our approach. Penalized algorithms  could be easily derived for other models such as the accelerated failure time model or the semiparametric transformation model for instance.

Although we have only presented asymptotic theoretical results, the simulation studies show clear evidence that our estimators outperform standard estimators for small sample sizes. Therefore, it would be of great interest to study their finite sample properties. However, such results involve deviation inequalities for non i.i.d. and non martingale empirical processes. %The main issue stands in studying the processes involving $M$ ($Z_n$ in the additive model and  in the Cox model) which, contrary to classical survival analysis, are not martingales. Deriving non asymptotic results of these processes is a real challenge and t
To our knowledge, no such results have yet been established in the context of recurrent events. 

 Another development of the present paper would be to establish results for the estimation of change-point locations and the number of change-points. Such results can be found for the change-point detection in the mean of a gaussian signal in~\cite{MR2796565}, for instance.
%Although we have only presented asymptotic theoretical results, it would be of great interest to study the finite sample properties of the estimators presented in this paper. As a matter of fact, the simulation studies show clear evidence that our estimators perform well for small sample sizes. However, such results are
%The results presented in this paper consider only asymptotic properties of our estimators. However, the simulation studies show clear evidence that our estimators perform well for small sample sizes. 
%Only asymptotic results were presented in this paper
%\begin{thebibliography}{7}

\appendix

\appendixone
\section*{Appendix: Proofs}\label{sec:Appendix}
Proofs of Lemma \ref{lem:relfond} to  \ref{lem} are in Supplementary Material.
\subsection*{A key relation}
%criteria of Equations~\eqref{eq:contrastMult} and~\eqref{eq:contrastAdd} are consequences of Lemma~A\ref{lem:relfond}.

%\texttt{ a verifier !!!!!!}
\begin{lemma}\label{lem:relfond}
Under Assumption~\ref{ass:independence}, for all $i=1,\ldots,n$
\begin{equation*}
\E \big( dN_i(t)\mid X_i(t), D_i \wedge C_i \geq  t,N_i(t)=s-1\big)=Y_i(t) \rho_0(t,s,X_i(t))dt.
\end{equation*}
\end{lemma}

\subsection*{Decomposition of the least squares criterion in the additive model}
The next proposition gives the details of the construction of the partial least squares in the additive model. One has to notice that the processes $Z_n(s)$ introduced below are centered which implies that finding a minimizer of $L_n^{PLS}$ is a natural way of estimating $\beta_0$ in model~\eqref{eq:modeladd}.
\begin{lemma}\label{prop:partialLS}
In the additive event-specific model~\eqref{eq:modeladd}, the partial least squares criterion~\eqref{eq:contrastAdd} can be rewritten as
\begin{align}\label{eq:LnAddProc}
L_n^{PLS}(\beta)&=
%  \frac 1n \sum_{i=1}^n \sum_{s=1}^B \int (X_i(t)\beta(s)-\bar X^s(t)\beta(s))^2
%  Y^s_i(t) \,dt \nonumber\\&\quad- \frac 2n  \sum_{i=1}^n \sum_{s=1}^B \int (X_i(t)\beta(s)-\bar X^s(t)\beta(s))(X_i(t)\beta_0(s)-\bar X^s(t)\beta_0(s))Y_i^s(t)dt\nonumber\\&\quad-
%\frac 2n  \sum_{i=1}^n \sum_{s=1}^B \int (X_i(t)\beta(s)-\bar X^s(t)\beta(s))  \ind{N_i(t)=s}dM_i(t)\nonumber\\&=
\sum_{s=1}^B\left\{ \beta(s)^\top \H_n(s) \beta(s)-2\beta(s)^\top \H_n(s) \beta_0(s)-2Z_n(s) \beta(s)\right\},
\end{align}
where
\[Z_n(s)=\frac 1n  \sum_{i=1}^n \sum_{s=1}^B \int \{X_i(t)-\bar X^s(t)\}\ind{N_i(t)=s}dM^s_i(t).\]
\end{lemma}
%The proof is given in Supplementary Material.
%\texttt{ cf. Lin wei yang ying p.724}

\subsection*{A technical lemma}
\begin{lemma}\label{lem} 
Let $\mathcal D[0,\tau]$ denotes the set of c\`adl\`ag functions on $[0,\tau]$ and let $F_n(\cdot)$ and $f(T,\delta,X(\cdot),N(\cdot ))$ be two random processes of bounded variation on $[0,\tau]$. Suppose that for all $z\in [0,\tau]$,
\[\mathbb E\left[\Big(\int_0^z f(T,\delta,X(t),N(t))dM^s(t)\Big)^{2}\right]<\infty.\]
We then have the following properties:\\
%%\begin{enumerate}[(i)]
%%\item \[\frac{1}{\sqrt n}\sum_{i=1}^n\int_0^zX_i(t)\mathds 1_{N_i(t)=s}dM_i(t)\]
%%\item If $\sup_{z\in [0,\tau]} F_n(z)=o_{\mathbb P}(1)$ then
%%\[\sup_{z\in [0,\tau]} \frac{1}{n}\sum_{i=1}^n\int_0^zF_n(t)dM(t)\underset{n\to \infty}{\overset{\mathbb P}{\longrightarrow}} 0.\]
%%\item 
%If $\sup_{t\in [0,\tau]} |F_n(t)-F(t)|=o_{\mathbb P}(1)$, where $F($\cdot$)$ is of bounded variation on $[0,\tau]$, then
%\[\frac{1}{\sqrt n}\sum_{i=1}^n\int_0^zF_n(t)f(T_i,\delta_i,X_i(t),N_i(t))dM_i(t)\]
%converges weakly in $\mathcal D[0,\tau]$ to a centered gaussian process with variance equal to
%\[\mathbb E\left[\Big(\int_0^zF(t)f(T,\delta,X(t),N(t))dM(t)\Big)^2\right].\]
%%%$\int_0^zF(t)\V[f(T,\delta,X(t),N(t))]dW(t)$ in $\mathcal D[0,\tau]$.\\
%
%\noindent OU\\
 
\begin{enumerate}[(i)]
\item If $f(T,\delta,X(\cdot),N(\cdot))$ is a random variable of bounded variation on $[0,\tau]$, then
\[\frac{1}{\sqrt n}\sum_{i=1}^n\int_0^zf(T_i,\delta_i,X_i(t),N_i(t))dM^s_i(t)\]
converges weakly in $\mathcal D[0,\tau]$ to a centered gaussian process with variance equal to
\[\mathbb E\left[\Big(\int_0^z f(T,\delta,X(t),N(t))dM^s(t)\Big)^{2}\right].\]
\item If $\sup_{t\in [0,\tau]} |F_n(t)-F(t)|=o_{\mathbb P}(1)$, where $F(\cdot)$ is a random process on $[0,\tau]$, then
\[\sup_{z\in [0,\tau]}\left\{\frac{1}{\sqrt n}\sum_{i=1}^n\int_0^z(F_n(t)-F(t))f(T_i,\delta_i,X_i(t),N_i(t))dM^s_i(t)\right\}=o_{\mathbb P}(1).\]
%\underset{n\to\infty}{\overset{\mathbb P}{\longrightarrow}} 0
\end{enumerate}
\end{lemma}
\subsection*{Proof of Theorem~\ref{theo:CVadd}}
\textsc{ Proof of 1.}  Let $\Gamma^{add}_n(\beta)$ be the quantity minimized by $\hat{\beta}_{\textsc{tv}/add}$ and introduce $\Gamma_{add}(\beta)= \sum_{s=1}^B \left[ \beta(s)^\top \H(s)\beta(s) -2 \bs h(s)\beta(s)  \right]$ where
\begin{align*}
%\H(s) & :=\int \mathbb E[Y^s(t)X(t)X(t)^\top]dt-\int\frac{(\mathbb E[Y^s(t)X(t)])^{\otimes 2}}{\mathbb E[Y^s(t)]}dt,\\
\bs h(s) & := \int \mathbb E\left[\ind{N(t)=s}X(t)dN(t)\right]-\int \frac{\mathbb E[Y^s(t)X(t)]}{\mathbb E[Y^s(t)]}\mathbb E[\ind{N(t)=s}dN(t)].
\end{align*}
Using Lemma~A\ref{lem:relfond} notice that $\bs h(s)=\beta_0(s)^\top\H(s)$ and consequently, $\argmin_{\beta}\Gamma_{add}=\beta_0$. Since the criterion to minimize is convex, the convergence in probability of $\hat{\beta}_{\textsc{tv}/add}$ to $\beta_0$ follows from the pointwise convergence of $\Gamma^{add}_n(\beta)$ towards $\Gamma_{add}(\beta)$. 
Now write:
\begin{align*}
%& \Big| L_n^{PLS}(\beta)+\frac{\lambda_n}{n}\sum_{j=1}^p \sum_{s=1}^B | \beta^j(s)-\beta^j(s-1)|-\Gamma(\beta)\Big|\\
& \Big|\Gamma^{add}_n(\beta) -\Gamma_{add}(\beta)\Big| \quad \leq  \Big| L_n^{PLS}(\beta)-\Gamma(\beta)\Big|+\frac{\lambda_n}{n}Bp\max_{s,j}| \beta^j(s)-\beta^j(s-1)|\\
& \quad \leq Bp^2\max_{j,k,s}|\beta^j(s)\beta^k(s)(\H_n^{j,k}(s)-\H^{j,k}(s))|+2Bp\max_{j,s} |\bs h_n^j(s)-\bs h^j(s)| |\beta^j(s)|+\frac{\lambda_n}{n}Bp
\end{align*}
and the result follows from the law of large number and the fact that $\lambda_n/n \to 0$ as n tends to infinity.

\textsc{ Proof of 2.} Define
\begin{equation*}
\Lambda^{add}_n(u)=\sum_{s=1}^B u(s)^\top \H_n(s) u(s)
-2 \sqrt{n}\sum_{s=1}^B Z_n(s)u(s)  +\lambda_n\sum_{j=1}^p\left(\TV(\beta^j_0+u^j/\sqrt{n})-\TV(\beta^j_0)\right)
\end{equation*}
%where $\TV(\beta)=\sum_{j=1}^p \sum_{s=2}^B | \beta^j(s)-\beta^j(s-1)|$ 
and notice that $\Lambda^{add}_n(u)$ is minimum at $u=\sqrt n(\hat \beta_{\textsc{tv}/add}-\beta_0)$.
Write
\begin{align*}
&\sqrt{n}\sum_{s=1}^B Z_n(s)u(s)  %&=\frac{1}{\sqrt n}\sum_{i=1}^n\int_0^{\tau}\sum_{s=1}^B u(s)^\top(X_i(t)-\bar X^s(t))\mathds 1_{N_i(t)=s}dM_i(t)
%\\& 
=\frac{1}{\sqrt n}\sum_{i=1}^n\int_0^{\tau}\sum_{s=1}^B \left(X_i(t)- \frac{\mathbb E[Y^s(t)X(t)]}{\mathbb E[Y^s(t)]}\right)u(s)\ind{N_i(t)=s}dM^s_i(t)
\\&
\quad -\frac{1}{\sqrt n}\sum_{i=1}^n\int_0^{\tau}\sum_{s=1}^B \left(\bar X^s(t) -\frac{\mathbb E[Y^s(t)X(t)]}{\mathbb E[Y^s(t)]}\right)u(s)\ind{N_i(t)=s}dM^s_i(t).
\end{align*}
Let $F_n(t)=\sum_s(\bar X^s(t)- \mathbb E[Y^s(t)X(t)]/\mathbb E[Y^s(t)])u(s)$ and $F(t)=0$. $F_n$ has bounded variation and from Lemma A\!\!\!~\ref{lem} (ii), the second term converges to $0$ in probability. Now, take $f(T_i,\delta_i,X_i(t),N_i(t))=\sum_s(X_i(t)- \mathbb E[Y^s(t)X(t)]/\mathbb E[Y^s(t)])u(s)\ind{N_i(t)=s}$ which is also a function of bounded variation. From Lemma A\!~\ref{lem} (i), the first term converges weakly towards a centered gaussian variable with variance equal to
\begin{align*}
& \mathbb E\left[\Big(\int_0^{\tau} \sum_{s=1}^B (X(t)- \mathbb E[Y^s(t)X(t)]/\mathbb E[Y^s(t)])u(s)\ind{N(t)=s}dM^s(t)\Big)^{2}\right]\\
&\quad = \sum_{s=1}^B u(s)^\top\mathbb E\left[\Big(\int_0^{\tau}(X(t)-\mathbb E[Y^s(t)X(t)]/\mathbb E[Y^s(t)])\ind{N(t)=s}dM^s(t)\Big)^{\otimes 2}\right]u(s).
\end{align*}
Then, note that $\sum_{s=1}^B u(s)^\top \H_n(s) u(s)$ converges to $\sum_{s=1}^B u(s)^\top \H(s) u(s)$, in probability 
%\[\sum_{s=1}^B u(s)^\top \H_n(s) u(s) \underset{n\to\infty}{\overset{\mathbb P}{\longrightarrow}} \sum_{s=1}^B u(s)^\top \H(s) u(s)\]
and $ \lambda_n\sum_j\left(\TV(\beta^j_0+u^j/\sqrt{n})-\TV(\beta^j_0)\right)/\lambda_0$ converges to
\begin{align*}
\sum_{j=1}^p \sum_{s=2}^B \left\{| \Delta u^j(s)|\ind{\Delta \beta^j(s)=0}\right.
%\\&
 %\qquad\quad 
 +\left.\text{sgn}(\Delta\beta_0^j(s))(\Delta u^j(s))\ind{\Delta \beta^j(s)\neq0}\right\}.
\end{align*}
Thus $\Lambda^{add}_n(u)$ converges to $\Lambda_{add}(u)$ in distribution. Since $\Lambda^{add}_n$ is convex and $\Lambda_{add}$ has a unique minimum, it follows that
$\sqrt n(\hat \beta_{\textsc{tv}/add} -\beta_0)$ converges to $\argmin_u \Lambda_{add}(u)$ in distribution.

\subsection*{Proof of Theorem~\ref{theo:CVmult}}
First define for $l=0,1$ or $2$ 
\[S_n^{(l)}(s,t,\beta)=\frac{1}{n}\sum_{i=1}^nY_i^s(t)X_i(t)^{\otimes l}\exp(X_i(t)\beta(s)).\]
% Following the arguments pages 305 and 306 of~\cite{FlemHar}, it can easily be shown that
% \[\sup_{t\in[0,\tau]}|S_n^{(l)}(s,t,\beta_0)-s^{(l)}(s,t,\beta_0)|\underset{n\to\infty}{\overset{\mathbb P}{\longrightarrow}} 0,\, \forall\,l=0,1,2,\]
% using the fact that the covariates process is of bounded variation (in particular, this assumption guarantees that $s^{(l)}(s,t,\beta_0)$ has a countable number of jumps).
Following the arguments in example VII.2.7 page 502 of~\cite{ABGK}, it can easily be shown that
\[\sup_{t\in[0,\tau]}|S_n^{(l)}(s,t,\beta_0)-s^{(l)}(s,t,\beta_0)|\underset{n\to\infty}{\overset{\mathbb P}{\longrightarrow}} 0,\, \forall\,l=0,1,2,\]
using the fact that the covariates process is of bounded variation (in particular, this assumption guarantees that $s^{(l)}(s,t,\beta_0)$ has a countable number of jumps).

\textsc{ Proof of 1.}  Let $\Gamma^{mult}_n(\beta)$ be the quantity minimized by $\hat{\beta}_{\textsc{tv}/mult}$ and introduce
%\[\Gamma^{mult}_n(\beta)= -\frac{1}{n}\sum_{s=1}^B \sum_{i=1}^n  \int  \left[X_i(t)^\top \beta(s)-\log (S_n^{(0)}(s,t,\beta))\right] Y_i^s(t)dN_i(t),\]
%%where
%%\[S_n^{(0)}(s,t,\beta)=\frac{1}{n}\sum_{i=1}^nY_i^s(t)\exp(X_i(t)^\top\beta(s))\]
%and notice that 
%\[\hat \beta_{\textsc{tv}/mult} \in \underset{\beta \in \mathbb R^{p\times B} } \argmin \left\{\Gamma^{mult}_n(\beta)+\frac{\lambda_n}{n}\sum_{j=1}^p \sum_{s=2}^B | \Delta\beta^j(s)|\right\}.\]
%%\beta^j(s)-\beta^j(s-1)
%%From classical results in Cox models (see for instance $\ldots$) it can be proved that
%%\[\sup_{t\in[0,\tau]}|S_n^{(0)}(s,t,\beta)-s^{(0)}(s,t,\beta)|\underset{n\to\infty}{\overset{\mathbb P}{\longrightarrow}} 0,\]
%%where 
%%\[s^{(0)}(s,t,\beta)=\mathbb E[Y^s(t)\exp(X(t)^\top \beta(s))].\]
%Introduce 
\begin{align*}
\Gamma_{mult}(\beta) & = -\sum_{s=1}^B  \int  \mathbb E\left[X(t) \beta(s)Y^s(t)dN(t)\right]+\sum_{s=1}^B\mathbb \int \log (s^{(0)}(s,t,\beta))\mathbb E\left[Y^s(t)dN(t)\right]\\
              & = -\sum_{s=1}^B  \int  \alpha_0(t,s) \left(\beta(s)^{\top}s^{(1)}(s,t,\beta_0)-\log (s^{(0)}(s,t,\beta))s^{(0)}(s,t,\beta_0)\right)dt,
\end{align*}
where the last equality follows from Lemma~A\ref{lem:relfond}. From similar arguments as in proof 1. of Theorem~\ref{theo:CVadd} and the uniform convergence with respect to $t$ of $S_n^{(0)}(s,t,\beta_0)$ towards $s^{(0)}(s,t,\beta_0)$, we get the pointwise convergence in probability of $\Gamma^{mult}_n(\beta)$ to $\Gamma_{mult}(\beta)$.
%\[|\Gamma^{mult}_n(\beta)-\Gamma_{mult}(\beta)|\underset{n\to\infty}{\overset{\mathbb P}{\longrightarrow}} 0.\]
Then, the consistency of $\hat \beta_{\textsc{tv}/mult}$ follows from the convexity of $\Gamma^{mult}_n(\beta)$ and the fact that $\argmin_{\beta} \,\Gamma_{mult}(\beta)=\beta_0$.%towards $\beta_0$ 

\textsc{ Proof of 2.} Consider the convex function
\[\Lambda^{mult}_n(u)=n\Gamma_n(\beta_0+u/\sqrt n)-n\Gamma_n(\beta_0)+\lambda_n\sum_{j=1}^p\left(\TV(\beta^j_0+u^j/\sqrt{n})-\TV(\beta^j_0)\right)\]
which is minimum at $u=\sqrt n(\hat \beta_{\textsc{tv}/mult}-\beta_0)$. Then from a Taylor expansion, one gets
\begin{align*}
&\Lambda^{mult}_n(u) =-\frac{\sqrt n}{n}\sum_{s=1}^B\sum_{i=1}^n \int \left(X_i(t)-\mathbf E_n(s,t,\beta_0)\right)Y_i^s(t)dN_i(t)u(s) \\
       & \quad +\frac{1}{2n}\sum_{s=1}^Bu(s)^\top \sum_{i=1}^n \int \mathbf V_n(s,t,\beta_0)Y_i^s(t)dN_i(t)u(s)%\\
       %&
       +\lambda_n\sum_{j=1}^p\left(\TV(\beta^j_0+u^j/\sqrt{n})-\TV(\beta^j_0)\right)+o_{\mathbb P}(1),%\quad 
\end{align*}
where
\begin{align*}
\mathbf E_n(s,t,\beta)  =\frac{S_n^{(1)}(s,t,\beta)}{S_n^{(0)}(s,t,\beta)},\;\;\;\;
%\mathbf V_n(s,t,\beta) & = \frac{S_n^{(2)}(s,t,\beta)}{S_n^{(0)}(s,t,\beta)}-\left(\frac{S_n^{(1)}(s,t,\beta)}{S_n^{(0)}(s,t,\beta)}\right)^{\otimes 2}.
\mathbf V_n(s,t,\beta)  = \frac{S_n^{(2)}(s,t,\beta)}{S_n^{(0)}(s,t,\beta)}-\mathbf E_n(s,t,\beta)^{\otimes 2}.
%S_n^{(1)}(s,t,\beta) & =\frac{1}{n}\sum_{i=1}^nY_i^s(t)X_i(t)\exp(X_i(t)^\top\beta(s)),\\
%S_n^{(2)}(s,t,\beta) & =\frac{1}{n}\sum_{i=1}^nY_i^s(t)X_i(t)\exp(X_i(t)^\top\beta(s))X_i(t)^\top.
\end{align*}
The uniform convergence with respect to $t$ of $S_n^{(0)}(s,t,\beta)$ and $S_n^{(2)}(s,t,\beta)$ towards $s^{(0)}(s,t,\beta_0)$ and $s^{(2)}(s,t,\beta_0)$ respectively and the law of large number give the convergence in probability of the term
\[\frac{1}{2n}\sum_{s=1}^Bu(s)^\top \sum_{i=1}^n \int \mathbf V_n(s,t,\beta_0)Y_i^s(t)dN_i(t)u(s)\]
towards
\[\frac{1}{2}\sum_{s=1}^Bu(s)^\top \int \mathbf v(s,t,\beta_0)\mathbb E[Y^s(t)dN(t)]u(s).\]
Notice that 
\[\sum_{i=1}^n\left(X_i(t)-\mathbf E_n(s,t,\beta_0)\right)Y_i^s(t)\alpha_0(t,s)\exp(X(t)\beta_0(t))dt=0\]
in order to rewrite the first term of $\Lambda^{mult}_n(u)$ as
\[-\frac{\sqrt n}{n}\sum_{s=1}^B\sum_{i=1}^n \int \left(X_i(t)-\mathbf E_n(s,t,\beta_0)\right)u(s)Y_i^s(t)dM^s_i(t).\]
From Lemma~\ref{lem}, the same kind of arguments as in the proof of Theorem~\ref{theo:CVadd} can be applied to conclude the proof.%obtain 

%the asymptotic weak convergence of this quantity towards a centered gaussian variable with variance equal to
%% \begin{align*}
%% & \mathbb E \left[\left(\sum_{s=1}^Bu(s)^\top \int \left(X(t)-\frac{s^{(1)}(s,t,\beta_0)}{s^{(0)}(s,t,\beta_0)}\right)Y^s(t)dM(t)\right)^2\right] \\
%%  & =\sum_{s=1}^B u(s)^\top \mathbb E \left[\left(\int \left(X(t)-\frac{s^{(1)}(s,t,\beta_0)}{s^{(0)}(s,t,\beta_0)}\right)Y^s(t)dM(t)\right)^{\otimes2}\right] u(s).
%% \end{align*}
%\begin{align*}
%%& \mathbb E \left[\left(\sum_{s=1}^Bu(s)^\top \int \left(X(t)-\mathbf e(s,t,\beta_0)\right)Y^s(t)dM(t)\right)^2\right] \\
 %%& =
 %\sum_{s=1}^B u(s)^\top \mathbb E \left[\left(\int \left(X(t)-\mathbf e(s,t,\beta_0)\right)Y^s(t)dM(t)\right)^{\otimes2}\right] u(s).
%\end{align*}
%We conclude as in the proof of Theorem~\ref{theo:CVadd}.

\bibliographystyle{biometrika}
\expandafter\ifx\csname natexlab\endcsname\relax\def\natexlab#1{#1}\fi
\bibliography{biblio}

\section*{Supplementary material}
\label{SM}
 Supplementary material %available at \Bka\ online
   includes a description of the algorithms, extended simulation study and additional analysis on the bladder tumour data of~\cite{byar}. It also contains proofs of Proposition~\ref{prop:partialLS} and Lemma~\ref{lem}.
%\pagebreak
%\include{Supp_Mat_recurrent_event_specific}
\end{document}